\documentclass[12pt]{article}
\usepackage{mathrsfs}
\usepackage{stmaryrd}
\usepackage{amsfonts,amsmath,amssymb,amscd}
\usepackage{shadow}
\usepackage{graphicx}
\usepackage{color}
\usepackage{longtable}
\usepackage{pstricks,multido}
\usepackage{hyperref}
\usepackage{qtree}
\allowdisplaybreaks

\newtheorem{thm}{Theorem}[section]

\newtheorem{conj}[thm]{Conjecture}
\def\qed{\hfill \rule{4pt}{7pt}}

\def\pf{\noindent {\it{Proof.} \hskip 2pt}}

\numberwithin{equation}{section}

\begin{document}
\begin{center}
{\large\bf
$s$-Inversion Sequences and $P$-Partitions of Type $B$}
 \end{center}

\begin{center}
{\small William Y.C. Chen$^1$, Alan J.X. Guo$^2$, Peter L. Guo$^3$
\\ Harry H.Y. Huang$^4$, Thomas Y.H. Liu$^5$}

\vskip 3mm
$^{1, 2, 3, 4, 5}$Center for Combinatorics, LPMC-TJKLC\\
Nankai University\\
Tianjin 300071, P.R. China
\\[3mm]

\vskip 4mm

$^1$chen@nankai.edu.cn, $^2$aalen@mail.nankai.edu.cn,
$^3$lguo@nankai.edu.cn\\
 $^4$hhuang@cfc.nankai.edu.cn,
$^5$lyh@cfc.nankai.edu.cn
\end{center}

\begin{abstract}
Given a sequence $s=(s_1,s_2,\ldots )$ of positive integers,
the inversion sequences with respect to $s$, or $s$-inversion sequences,
 were introduced by Savage and Schuster
in their study of lecture hall polytopes.  A sequence $(e_1,e_2,\ldots,e_n)$ of nonnegative
integers is called  an $s$-inversion sequence of length $n$
 if $0\leq e_i < s_i$ for $1\leq i\leq n$.
 Let $I(n)$ be the set of $s$-inversion sequences of length $n$
 for $s=(1,4,3,8,5,12,\ldots)$,
that is,  $s_{2i}=4i$
and $s_{2i-1}=2i-1$ for  $i\geq1$, and let $P_n$ be the set of signed permutations on $\{1^2,2^2,\ldots,n^2\}$.
Savage and  Visontai conjectured  that when $n=2k$,
the ascent number over  $I_n$ is equidistributed with the  descent number over $P_k$.  For a  positive integer $n$, we use type $B$
$P$-partitions to give a characterization of  signed  permutations
over which the descent number is equidistributed with the ascent number over  $I_n$. When $n$ is even, this confirms the
conjecture of Savage and Visontai.
Moreover, let $I'_n$ be the set of $s$-inversion sequences of length $n$ for $s=(2,2,6,4,10,6,\ldots)$,
that is,  $s_{2i}=2i$ and $s_{2i-1}=4i-2$ for $i\geq1$. We find
a set of signed   permutations
over which the descent number is equidistributed
with the ascent number over $I'_n$.
\end{abstract}



\section{Introduction}

 The notion of $s$-inversion sequences
was introduced by Savage and Schuster  \cite{SavageM}
in their study of lecture hall polytopes.
Let $s=(s_1,s_2,\ldots)$ be a sequence of positive integers.
An  inversion sequence of length $n$ with respect to $s$, or
an $s$-inversion sequence of length $n$, is a sequence
$e=(e_1, e_2, \ldots, e_n)$ of nonnegative integers
 such that $0\leq e_i<s_i$ for $1\leq i\leq n$.
An ascent of an $s$-inversion sequence $e=(e_1, e_2, \ldots, e_n)$ is defined
to be an integer $i\in \{0,1,\ldots,n-1\}$ such that
\[\frac{e_i}{s_i}<\frac{e_{i+1}}{s_{i+1}},\]
under the assumption that  $e_0=0$ and $s_0=1$.
The ascent number $\mathrm{asc}(e)$ of $e$ is meant to be the
number of ascents of $e$.

The generating function of the ascent number over $s$-inversion
sequences can be viewed as a generalization of the Eulerian
polynomial for permutations, since the ascent number over the $s$-inversion sequences of length $n$ for $s=(1,2,3,\ldots)$
is equidistributed with the descent number over the permutations on $\{1,2,\ldots,n\}$, see Savage and Schuster \cite{SavageM}.
For an inversion sequence $e=(e_1, e_2, \ldots, e_n)$  with respect to $s=(s_1,s_2,\ldots,)$,  let
\[
\mathrm{amaj}(e)=\sum_{i\in \mathrm{Asc}(e)}(n-i),
\]
and
\[
\mathrm{lhp}(e)=-|e|+\sum_{i\in \mathrm{Asc}(e)}(s_{i+1}+\cdots+s_n),
\]
where $\mathrm{Asc}(e)$ is the set of ascents of $e$,
and $|e|=e_1+e_2+\cdots+e_n$ is the weight of $e$.
Savage and Schuster \cite{SavageM} showed that
the multivariate  generating function for  the ascent number $\mathrm{ase}(e)$, the major index $\mathrm{amaj}(e)$, the lecture
hall statistic $\mathrm{lhp}(e)$ and the weight $|e|$ is related to
 the  Ehrhart series of  $s$-lecture hall
polytopes and the  generating function  of
$s$-lecture hall partitions.

Savage and Visontai \cite{Savage} found a connection
between the generating function  of the ascent number
over $s$-inversion sequences of length $n$  and a conjecture
of Brenti \cite{Brenti} on the real-rootedness of
Eulerian polynomials of finite Coxeter groups.
The real-rootedness of the Eulerian polynomial of type $A$
was known to Frobenius \cite{Frobenius}, see also \cite{Bona,Haglund}. Brenti \cite{Brenti} proved the
real-rootedness of the
Eularian polynomials of Coxeter groups of  type $B$ and exceptional Coxeter groups.
For the sequence $s=(2,4,6,\ldots)$ and an $s$-inversion sequence $e=(e_1,e_2,\ldots,e_n)$, Savage and Visontai \cite{Savage} defined the type $D$ ascent set of $e$ as given by
\[\mathrm{Asc}_D(e)=\left\{\,i\,\Big|\,
\frac{e_i}{i}<\frac{e_{i+1}}{i+1},\, 1\leq i\leq n-1\right\}\cup \{0\,|\,\mathrm{if}\, 2e_1+e_2\geq 3\}.\]
Let $T_n(x)$ be
the generating function of the type $D$ ascent number over $s$-inversion sequence of length $n$ for $s=(2,4,6,\ldots)$.
For example, $T_3(x)=2(x^3+11x^2+11x+1)$.  Let $D_n(x)$ be the $n$-th Eulerian polynomial of type $D$. Recall that the type $D$ Coxeter group of rank $n$, denoted $D_n$,  is the group of even-signed permutations
on $\{1,2,\ldots,n\}$, see Bj\"orner and Brenti \cite{BB}. The  descent set $\mathrm{Des}_D(\sigma)$ of an even-signed
permutation $\sigma=\sigma_1\sigma_2\cdots \sigma_n\in D_n$ is
defined by
\[\mathrm{Des}_D(\sigma)=\left\{\,i\,|\,
\sigma_i>\sigma_{i+1}, 1\leq i\leq n-1\right\}\cup \{0\,|\,\mathrm{if}\, \sigma_1+\sigma_2<0\} .\]
The descent number of $\sigma$ is meant to be the number of
descents in $\mathrm{Des}_D(\sigma)$.
Let $D_n(x)$ denote the generating function of the descent number over $D_n$.
Savage and Visontai \cite{Savage}
showed that $T_n(x)=2D_n(x)$.
 By proving that  $T_n(x)$ has only real roots for $n\geq 1$, they deduced the real-rootedness of
the Eulerian polynomials of type $D$ and
settled the last unsolved case of the conjecture of Brenti \cite{Brenti}.

Savage and Visontai \cite{Savage} proved  that for
 any sequence $s$ of positive integers and any positive integer $n$,
the generating
function of the ascent number over  $s$-inversion sequences
of length $n$ has
only real roots. Let  $I_n$ denote the set of  $s$-inversion sequences of length $n$ for the specific sequence $s=(1,4,3,8,5,12,\ldots)$, that is, for $i\geq1$, $s_{2i}=4i$
 and $s_{2i-1}=2i-1$.  Let $P_n$ denote the set of
signed permutations on the multiset
$\{1^2,2^2,\ldots,n^2\}$.
Savage and Visontai \cite{Savage} posed the following conjecture, which
implies the real-rootedness of
the generating function of the descent number over $P_n$.

\begin{conj}[\mdseries{\cite[Conjecture 3.27]{Savage}}]\label{thatConj}
For $n\geq 1$, the descent number over $P_n$
is equidistributed with the ascent number over $I_{2n}$.
\end{conj}

In this paper, we give  a proof  of
Conjecture \ref{thatConj}. Let $P_n(x)$ denote the generating
function of the descent number over the set $P_n$ of signed permutations on $\{1^2,2^2,\ldots,n^2\}$, and let $I_{n}(x)$ denote the generating function of  the ascent number
over $I_{n}$. Savage and Schuster \cite{SavageM} found a relation
for $I_{n}(x)$.
We show that the generating function $P_n(x)$
equals  the generating function of the descent number over  linear extensions of
certain signed labeled forests.
By using
$P$-partitions of type $B$  introduced by
Chow \cite{Chow}, we show  that the generating function
 for the descent number over linear
 extensions satisfies the same
relation as $I_{2n}(x)$. Thus the generating function $P_n(x)$ satisfies the same
relation as $I_{2n}(x)$.
This proves Conjecture \ref{thatConj}.

We also find  characterizations of signed permutations such that the
descent number is equidistributed with the ascent number over three other classes of $s$-inversion sequences. To be specific, we show that the descent number over the set of  signed permutations on the  multiset $\{1^2,2^2,\ldots,(n-1)^2,n\}$ such that  $n$ is assigned a minus sign
 is equidistributed with the ascent number over $I_{2n-1}$.
For
$s=(2,2,6,4,10,6,\ldots)$, that is, for $i\geq1$, $s_{2i}=2i$ and $s_{2i-1}=4i-2$, let  $I_n'$ denote the set of  $s$-inversion sequences of length $n$. We show that the descent number over $P_n$ is
equidistributed with the ascent number over $I_{2n}'$ and
the descent number over the set of  signed permutations on
 $\{1^2,2^2,\ldots,(n-1)^2,n\}$ is
equidistributed with the ascent number over $I_{2n-1}'$.

\section{Proof of Conjecture \ref{thatConj}}\label{Ppartition}

In this section, we  present a proof of Conjecture \ref{thatConj}
by establishing a connection between the generating function $P_n(x)$
of the descent number over $P_n$
and the generating function of the descent number over linear extensions of certain signed  labeled forests with $2n$ vertices.  Let $F_{n}$ be the plane forest with $n$
trees containing exactly two vertices, and let $F_n(x)$
denote  the generating
function  of the descent number over linear extensions of $(F_n,w)$, where $w$
ranges over certain signed labelings of $F_n$.
Keep in mind that a plane forest means a set of
plane trees that are arranged in linear order.
We shall show that $F_n(x)= P_n(x)$. On the other hand,
by using  $P$-partitions of type $B$ introduced by Chow
\cite{Chow}, we obtain a relation for $F_n(x)$.
Savage and Schuster \cite{SavageM} have shown that  the same relation is   satisfied by  the generating function
$I_{2n}(x)$, so we get $I_{2n}(x)=F_n(x)$.
This confirms Conjecture \ref{thatConj}, that is, $P_n(x)=I_{2n}(x)$.

\begin{figure}[h,t]
\setlength{\unitlength}{0.5mm}
\begin{center}
\begin{picture}(100,40)
\put(0,15){\circle*{2}}\put(10,30){\circle*{2}}\put(20,15){\circle*{2}}
\put(0,15){\line(2,3){10}}\put(10,30){\line(2,-3){10}}
\put(50,15){\circle*{2}}\put(50,30){\circle*{2}}
\put(50,30){\line(0,-1){15}}
\put(80,0){\circle*{2}}\put(90,30){\circle*{2}}\put(90,15){\circle*{2}}
\put(100,0){\circle*{2}}\put(80,0){\line(2,3){10}}\put(90,15){\line(2,-3){10}}
\put(90,30){\line(0,-1){15}}
\put(-6,5){$-1$}\put(18,5){$7$}\put(3,30){$4$}
\put(47,5){$2$}\put(38,30){$-9$}
\put(68,-5){$-5$}\put(102,-5){$6$}\put(93,12){$-8$}\put(93,29){$3$}
\end{picture}
\end{center}\caption{A signed labeled forest on $\{1,2,\ldots,9\}$}\label{clg}
\end{figure}
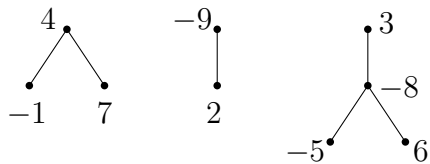

Let us give an overview of linear extensions of  a signed labeled forest.
Let $F$ be a plane forest with $n$ vertices,
and let $S$ be a set of $n$ distinct positive integers.
A labeling of $F$ on $S$ is an assignment of the elements in $S$
to the vertices of $F$ such that each element in $S$ is assigned
to only one vertex. A signed labeling of $F$ on $S$ is
a labeling of $F$ on $S$ with each label possibly associated with  a  minus sign.
For example, Figure \ref{clg} illustrates a signed labeled forest on $\{1,2, \ldots, 9\}$.
We  use $(F,w)$ to stand for a plane forest $F$ associated with
a signed labeling $w$.

Linear extensions of
a signed labeled forest $(F,w)$ are defined
 based on linear extensions of  $F$.
For a plane forest $F$ with $n$ vertices, say $x_1,x_2,\ldots,x_n$,
a linear extension of $F$ is a permutation
$x_{i_1}x_{i_2}\cdots x_{i_n}$ of the vertices of $F$
 such that   $x_{i_j}<_Fx_{i_k}$ implies $j<k$, where $<_F$ is the order relation of $F$. Let $\mathcal{L}(F)$ denote the
set of linear extensions of $F$.
Then
the set of linear extensions of $(F,w)$ is defined as
\[\mathcal{L}(F,w)=\left\{w(x_{i_1})w(x_{i_2})\cdots w(x_{i_n})\,|\,
x_{i_1}x_{i_2}\cdots x_{i_n}\in \mathcal{L}(F)\right\},\]
where $w(x)$  denotes the label of a vertex $x$ of $F$.

Notice that a linear extension of $(F,w)$
is a signed permutation on the labeling set $S$ of $F$.
We next define the generating function $F_n(x)$ of the descent
number over linear extensions of the plane forest $F_n$
associated with certain signed labelings.
Let us recall  the descent number of a signed permutation on a multiset.
A signed permutation on a multiset $M$ is a permutation on $M$
for which each element is possibly assigned a minus sign.
For example, $\overline{3}1\overline{2}13$ is a signed permutation on $\{1^2,2,3^2\}$, where we use
a bar to indicate that an element is assigned a minus sign.
The descent set of a signed permutation $\sigma=\sigma_1\sigma_2\cdots \sigma_n$
is defined as
\begin{equation}\label{descent-2}
\{i\,|\, \sigma_i>\sigma_{i+1},\, 1\leq i\leq n-1\}\cup \{0\,|\, \text{if $\sigma_1<0$}\},
\end{equation}
see Savage and Visontai \cite{Savage}.
However, for  the purpose of this paper, we  use the
following reformation of the descent set of a signed permutation $\sigma=\sigma_1\sigma_2\cdots \sigma_n$
\begin{equation}\label{descent}
\mathrm{Des}_{B}(\sigma)=\{i\,|\, \sigma_i>\sigma_{i+1},\,1\leq i\leq n-1\}\cup \{n\,|\, \text{if $\sigma_n>0$}\}.
\end{equation}
The number of descents in $\mathrm{Des}_{B}(\sigma)$
 is referred to as the descent number  of $\sigma$, denoted $\mathrm{des}_{B}(\sigma)$.
In fact, via the    bijection
\[\sigma=\sigma_1\sigma_2\cdots \sigma_n\longmapsto \sigma'=(-\sigma_n)(-\sigma_{n-1})\cdots (-\sigma_1),\]
we see that the descent numbers defined by \eqref{descent-2} and $\mathrm{Des}_{B}(\sigma)$  are equidistributed over  signed permutations on any   multiset.

To  define the generating function $F_n(x)$, we introduce some  specific singed labelings for the plane forest $F_n$.
We write $T_1,T_2,\ldots, T_n$ for the $n$ trees of $F_n$, which are listed from left to right.
For $1\leq i\leq n$, let $u_i$  denote the root of $T_i$
 and let  $v_i$ denote the child of $u_i$.
Let $w$ be a signed labeling of $F_n$ on $\{1,2,\ldots,2n\}$
such that for $1\leq i\leq n$,
\[\left\{|w(u_i)|,\,|w(v_i)|\right\}=\{2i-1,2i\}.\]
Let $w_i$ $(1\leq i\leq n)$ be the signed  labeling of
$T_i$ induced by $w$.
There are eight possibilities for each $w_i$. However,
for the purpose of establishing the following equidistribution theorem, we need only four cases as given below:
\begin{itemize}
\item[] Case 1: $w_i(u_i)=2i$ and $w_i(v_i)=2i-1$;
\item[] Case 2: $w_i(u_i)=2i$ and $w_i(v_i)=\overline{2i-1}$;
\item [] Case 3: $w_i(u_i)=\overline{2i}$ and $w_i(v_i)=2i-1$;
\item [] Case 4: $w_i(u_i)=\overline{2i-1}$ and $w_i(v_i)=\overline{2i}$.
\end{itemize}
For the $j$-th case, we say that $w_i$ is of type $j$.
Let $L(F_n)$ denote the set of  signed labelings  of $F_n$ such that the induced labeling $w_i$ of $T_i$ is one of the above four types.
We shall show that the descent number over $P_n$ is equidistributed  with the descent number over the set of  linear extensions of $(F_n,w)$, where $w$ ranges over the set $L(F_n)$.
Specifically, define
\[F_n(x)=
\sum_{w\in L(F_n)}\ \sum_{\sigma\in\mathcal{L}(F_n,w)}x^{\mathrm{des}_B(\sigma)}.
\]
We have the following equidistribution theorem.

\begin{thm}\label{al}
For   $n \geq 1$, we have
\[F_n(x)=P_n(x).\]
\end{thm}

\pf We proceed to construct a descent preserving bijection $\phi$ from
the set
\[
\left\{\sigma\in \mathcal{L}(F_n,w)\,|\, w\in L(F_n)\right\}
\]
to the set $P_n$ of signed permutations on $\{1^2,2^2,\ldots,n^2\}$.
Let $\sigma=\sigma_1\sigma_2\cdots \sigma_{2n}$ be a linear extension  in
$\mathcal{L}(F_n,w)$, where $w\in L(F_n)$.
Define $\phi(\sigma)=\tau=\tau_1\tau_2\cdots \tau_{2n}$  as follows.
For $1\leq j\leq 2n$, $\tau_j$ has the same sign as $\sigma_j$, and
$|\tau_j|=i$   if $|\sigma_j|=2i-1$
or   $|\sigma_j|=2i$.
It is routine to check   that
$\phi$ is a bijection.
Moreover, it is readily verified  that
$j\in \{1,2,\ldots,2n\}$ is a descent
of $\sigma$ if and only if it is a descent
of $\tau$.
This completes the proof.
\qed

The next theorem gives an expression for the
generating function $F_n(x)$.

\begin{thm}\label{new-2}
For $n\geq 1$, we have
\begin{equation}\label{new-1}
\frac{F_n(x)}{(1-x)^{2n+1}}
    =\sum\limits_{t\geq0}(t+1)^n(2t+1)^nx^t.
\end{equation}
\end{thm}

To prove  Theorem \ref{new-2}, we need  a decomposition of $P$-partitions of type $B$ into $\sigma$-compatible maps due to Chow \cite{Chow}, where $\sigma$ is a linear extension of $P$. When the poset $P$ is associated with an ordinary labeling,  type $B$ $P$-partitions reduce to   ordinary $P$-partitions introduced by Stanley \cite{Stanley-97}.
To make a connection to Theorem \ref{new-2},
 it is sufficient to consider the case when $P$ is a plane forest.
In this case, we do not need the structure of
$P$-partitions of type $B$ in
full generality as given by Chow. For the case when
   $P$ is a plane forest, a $P$-partition of type $B$
was described by Chen, Gao and Guo \cite{chen}.

Let $F$ be a plane forest, and
$w$ be a signed labeling of $F$. Let $\mathbb{N}$
be the set of nonnegative integers.
 A  $(F,w)$-partition of type $B$ is a map $f$ from the set of vertices of $ F$ to
$\mathbb{N}$ that satisfies the following conditions:
\begin{itemize}
\item[(1)] $f(x)\leq f(y)$ if $x\geq_F y$;
\item[(2)] $f(x)<f(y)$ if $x>_F y$ and $w(x)<w(y)$;
\item[(3)] $f(x)\geq 1$ if $x$ is a root of $F$ with $w(x)>0$.
\end{itemize}

Analogous to the decomposition of ordinary $P$-partitions
given by Stanley \cite{Stanley-97},
Chow \cite{Chow} showed that type $B$ $(F,w)$-partitions
can be decomposed into $\sigma$-compatible maps, where $\sigma$
is a linear extension of $(F,w)$.
For a signed permutation $\sigma=\sigma_1\sigma_2\cdots\sigma_n$,
 a $\sigma$-compatible map $g$ is a map  from $\{\sigma_1,\sigma_2,\ldots,\sigma_n\}$ to $\mathbb{N}$ that satisfies the following conditions:
\begin{itemize}
\item[(1)] $g(\sigma_1)\geq g(\sigma_2)\geq\cdots\geq g(\sigma_n)$;
\item[(2)] For $i\in \{1,2,\ldots,n-1\}$, $g(\sigma_i)>g(\sigma_{i+1})$ if $\sigma_i>\sigma_{i+1}$;
\item[(3)] $g(\sigma_n)\geq 1$ if $\sigma_n>0$.
\end{itemize}
Let  $A(F,w)$ denote  the set of type $B$ $(F,w)$-partitions,
and let $A_{\sigma}$ denote the set of   $\sigma$-compatible maps.
The following decomposition is due to Chow \cite{Chow}.

\begin{thm}[\mdseries{\cite[Theorem 2.1.4]{Chow}}]
Let $F$ be a plane forest associated with a signed labeling $w$.
Then
\begin{equation}\label{dec1}
A(F,w)=\bigcup_{\sigma\in
\mathcal{L}(F,w)}\,A_\sigma.
\end{equation}
\end{thm}

For a nonnegative integer $t$, let $\Omega_{F}(w,t)$
denote the number of  type $B$
$(F,w)$-partitions $f$  such that
$f(x)\leq t$ for any $x\in F$.
When $w$ is an ordinary labeling,
Stanley \cite{Stanley-97} has established a relation between the generating function of the descent number over linear extensions of $(F,w)$ and the generating function of $\Omega_{F}(w,t)$. For  signed labeled forests, we have the following  relation.

\begin{thm}\label{pppt}
Let $F$ be a plane forest with $n$ vertices, and $w$ be a signed labeling of $F$ on
$\{1,2,\ldots,n\}$. Then
\begin{equation}\label{eqpartitionanddes2}
\frac{\sum_{\sigma\in\mathcal{L}(F,w)}x^{\mathrm{des}_B(\sigma)}}{(1-x)^{n+1}}
=\sum\limits_{t\geq0}\Omega_F(w,t)\,x^t.
\end{equation}
\end{thm}

\pf  We essentially follow the proof of Stanley \cite{Stanley-97} for  ordinary $P$-partitions.
For a signed permutation $\sigma=\sigma_1\sigma_2\cdots\sigma_n$ on $\{1,2,\ldots,n\}$, let $\Omega_\sigma(t)$ denote  the number of $\sigma$-compatible maps
$g$ with $g(\sigma_1)\leq t$. For any
 linear extension $\sigma$ in $\mathcal{L}(F,w)$, in view of the decomposition \eqref{dec1},
 relation \eqref{eqpartitionanddes2} can be deduced from
   the following relation
\begin{equation}\label{ppp}
\sum\limits_{t\geq0}\Omega_\sigma(t)\,x^t=
\frac{x^{\mathrm{des}_B(\sigma)}}{(1-x)^{n+1}}.
\end{equation}

For $1\leq i\leq n$, let
$d_i$ denote the number of descents of $\sigma$ that are greater
than or equal to $i$, that is,
\[d_i=|\{j\,|\,\sigma_j\geq \sigma_{j+1},\,i\leq j\leq n-1\}\cup \{n\,|\, \text{if $\sigma_n>0$}\}|.\]
Setting $\lambda_i=g(\sigma_i)-d_i$, we are led to a one-to-one correspondence
between the set of $\sigma$-compatible maps
$g$ with $g(\sigma_1)\leq t$ and the set of  partitions $(\lambda_1,\lambda_2,\ldots,\lambda_n)$ with $\lambda_1\leq t-\mathrm{des}_B(\sigma)$,
where the latter is    counted by
\[{n+t-\mathrm{des}_B(\sigma)\choose n},\]
see,  for example,  Stanley \cite{Stanley-97}.
Thus,
\begin{align*}
\sum\limits_{t\geq0}\Omega_\sigma(t)\,x^t&=\sum\limits_{t\geq0}{n+t-\mathrm{des}_B(\sigma)\choose n}x^t\\[5pt]
&=\frac{x^{\mathrm{des}_B(\sigma)}}{(1-x)^{n+1}},
\end{align*}
which agrees with \eqref{ppp}. This completes the proof.
\qed

We are now ready to prove Theorem \ref{new-2}.

\noindent
\textit{Proof of Theorem \ref{new-2}.}
By Theorem \ref{pppt},  the assertion  \eqref{new-1}
is equivalent to the following relation
\begin{equation}\label{al-1}
\sum_{w\in L(F_n)}\Omega_{F_n}(w,t)=\left((t+1)(2t+1)\right)^n.
\end{equation}
Let $T_1,T_2,\ldots,T_n$ be the $n$ trees of $F_n$ listed from left to right. Thus, for any signed labeled forest $(F_n,w)$, we have
\begin{equation}\label{lgt}
\Omega_{F_n}(w,t)=\prod _{i=1}^n \Omega_{T_i}(w_i,t),
\end{equation}
where $w_i$ is the signed labeling of $T_i$ induced by $w$.
Recall that for  a signed labeling $w$ in $L(F_n)$, each induced labeling $w_i$ has four choices.
For $1\leq j\leq 4$, let $w_{i}^{(j)}$ be the signed labeling
of $T_i$ that is of type $j$, so that the left-hand side of \eqref{al-1} can be rewritten as
\begin{equation}\label{ohn}
\sum_{w\in L(F_n)}\Omega_{F_n}(w,t)=\prod _{i=1}^n \left(\Omega_{T_i}(w_{i}^{(1)},t)+\Omega_{T_i}(w_{i}^{(2)},t)+
\Omega_{T_i}(w_{i}^{(3)},t)+\Omega_{T_i}(w_{i}^{(4)},t)\right).
\end{equation}
We claim that for any $1\leq i \leq n$,
\begin{equation}\label{skn}
\Omega_{T_i}(w_{i}^{(1)},t)+\Omega_{T_i}(w_{i}^{(2)},t)+
\Omega_{T_i}(w_{i}^{(3)},t)+\Omega_{T_i}(w_{i}^{(4)},t)=(t+1)(2t+1).
\end{equation}
Assume that $u_i$ is the root of $T_i$ and $v_i$ is the child of $u_i$.
Let $f$ be a type $B$ $(T_i,w_{i}^{(j)})$-partition such that $f(v_i)\leq t$.
Then we have
\begin{align*}
&0<f(u_i)\leq f(v_i)\leq t, \ \ \ \ \text{if  $j=1$};\\[3pt]
&0<f(u_i)\leq f(v_i)\leq t, \ \ \ \ \text{if  $j=2$};\\[3pt]
&0\leq f(u_i)<f(v_i)\leq t, \ \ \ \ \text{if  $j=3$};\\[3pt]
&0\leq f(u_i)\leq f(v_i)\leq t, \ \ \ \ \text{if $j=4$}.
\end{align*}
It follows that
\begin{equation}\label{331}
\Omega_{T_i}(w_{i}^{(1)},t)=\Omega_{T_i}(w_{i}^{(2)},t)=
\Omega_{T_i}(w_{i}^{(3)},t)={t+1 \choose 2}
\end{equation}
and
\begin{equation}\label{332}
\Omega_{T_i}(w_{i}^{(4)},t)={t+2 \choose 2}.
\end{equation}
Hence we obtain  \eqref{skn}, completing the proof.
\qed

In addition to Theorem  \ref{al} and  Theorem \ref{new-2},
a formula of  Savage and Schuster \cite{SavageM} is   needed to prove Conjecture \ref{thatConj}.
Recall that  $I_n$  is the set of  $s$-inversion sequences of length $n$ for the specific sequence $s=(1,4,3,8,5,12,\ldots)$, and $I_n(x)$
is the generating function of the ascent number over $I_n$.
Savage and Schuster \mdseries{\cite[Theorem 13]{SavageM}} showed that for $n\geq1$,
\begin{equation}\label{mainequation1}
\frac{I_n(x)}{(1-x)^{n+1}}=
\sum\limits_{t\geq0}(t+1)^{\lceil\frac{n}{2}\rceil}(2t+1)^{\lfloor\frac{n}{2}\rfloor}x^t.
\end{equation}
Replacing $n$ by $2n$ in \eqref{mainequation1}, we get
\begin{equation}\label{bkg}
\frac{I_{2n}(x)}{(1-x)^{2n+1}}
    =\sum\limits_{t\geq0}(t+1)^n(2t+1)^nx^t.
\end{equation}
Comparing \eqref{bkg} with \eqref{new-1}, we obtain $F_n(x)=I_{2n}(x)$.  By  Theorem \ref{al} we arrive at $P_n(x)=I_{2n}(x)$, completing the proof of Conjecture \ref{thatConj}.

\section{Further equidistributions}

In this section, we give characterizations of three sets
of signed permutations over which the descent number
is equidistributed with the ascent number over the sets  $I_{2n-1}$, $I'_{2n-1}$ and $I_{2n}'$ respectively.

Recall that $I_{2n-1}$ stands for the set of  $s$-inversion sequences of length $2n-1$ for $s=(1,4,3,8,5,12,\ldots)$, and
 $I_n'$ stands for the set of  $s$-inversion sequences of length $n$ for $s=(2,2,6,4,10,6,\ldots)$.
Let $U_n$ be the set of  signed permutations on
the multiset $\{1^2,2^2,\ldots,(n-1)^2,n\}$, and
let $V_n$ be the subset of $U_n$ consisting of
signed permutations such that the element $n$ carries a minus sign. We show that the descent number over $V_n$
is equidistributed with the ascent number over $I_{2n-1}$, the descent number over $U_n$
is equidistributed with the ascent number over $I_{2n-1}'$, and
the descent number over $P_n$
is equidistributed with the ascent number over $I_{2n}'$.

\begin{thm}\label{theorem-2}
For $n\geq 1$, we have
\begin{equation}\label{main-2}
  \sum\limits_{\sigma\in V_n}x^{\mathrm{des}_{B} (\sigma)}
  =\sum\limits_{e\in I_{2n-1}}x^{\mathrm{asc}(e)}.
\end{equation}
\end{thm}

\pf
Let
 $V_n(x)$ denote the sum on the left-hand side of \eqref{main-2}.
 Since a relation on the generating function $I_{2n-1}(x)$ is given by \eqref{mainequation1},
 it suffices to show that $V_n(x)$ satisfies the same
 relation as $I_{2n-1}(x)$, that is, for $n\geq 1$,
\begin{equation}\label{main-equation4}
\frac{V_n(x)}{(1-x)^{2n}}
=\sum\limits_{t\geq 0}(t+1)^n(2t+1)^{n-1}x^t.
\end{equation}

 We claim that $V_n(x)$ equals  the generating
function of the descent number over linear extensions of certain signed labeled
forests. To this end, let $F_n'$ be the plane forest which is obtained from  $F_{n-1}$ by adding
a single vertex as the rightmost component.
Let $T_1,T_2,\ldots,T_{n-1}$ denote  $n-1$ trees of $F_{n-1}$,
and let $T_n$ denote a single vertex.
Hence  $F_n'$ consists of plane trees $T_1, \ldots, T_{n-1}, T_n$.  Write  $L(F_n')$ for the set of  signed labelings $w$ of $F_n'$ such that $w(T_n)=-{(2n-1)}$ and the induced signed labeling of $w$ on $F_{n-1}$ belongs to $L(F_{n-1})$.
Set
\begin{equation}\label{nsk}
F_n'(x)=\sum_{w\in L(F_n')}\,\sum_{\sigma\in\mathcal{L}(F_n',w)}x^{\mathrm{des}_B(\sigma)}.
\end{equation}
Using the same reasoning as in the proof of Theorem \ref{al}, one can construct a descent preserving bijection between the set
\[
\left\{\sigma\in \mathcal{L}(F_n',w)\,|\, w\in L(F_n')\right\}
\]
and the set $V_n$. Hence we get
\begin{equation}\label{kll}
F_n'(x)=V_n(x),
\end{equation}
so that \eqref{main-equation4} can be rewritten as
\begin{equation}\label{eq-4}
\frac{F_n'(x)}{(1-x)^{2n}}
=\sum\limits_{t\geq 0}(t+1)^n(2t+1)^{n-1}x^t.
\end{equation}
Applying Theorem \ref{pppt} to  the set of
 signed labeled forests $(F_n',w)$ with $w\in L(F_n')$,
we obtain
\[
\frac{F_n'(x)}{(1-x)^{2n}}=\sum\limits_{t\geq 0}
\sum_{w\in L(F_n')}\Omega_{F_n'}(w,t)\,
x^t.
\]
Thus,  \eqref{eq-4} can be deduced from the
following relation
\begin{equation}\label{al-3}
\sum_{w\in L(F_n')}\Omega_{F_n'}(w,t)=(t+1)^n(2t+1)^{n-1}.
\end{equation}

Notice  that for  a signed labeling $w$ in $L(F_n')$, each induced labeling $w_i$ of $T_i$ for $1\leq i\leq n-1$  has four types, and the induced labeling $w_n$ of $T_n$ satisfies $w_n(T_n)=-(2n-1)$. For $1\leq i\leq n-1$ and $1\leq j\leq 4$,
let $w_{i}^{(j)}$ be the signed labeling
of $T_i$ that is of type $j$.
In the proof of Theorem \ref{new-2},
it has been shown that for $1\leq i\leq n-1$,
\begin{equation}\label{agag}
 \Omega_{T_i}(w_{i}^{(1)},t)+\Omega_{T_i}(w_{i}^{(2)},t)+
\Omega_{T_i}(w_{i}^{(3)},t)+\Omega_{T_i}(w_{i}^{(4)},t)=(t+1)(2t+1).
\end{equation}
Moreover, it is clear that $\Omega_{T_n}(w_n,t)=t+1$.
Hence we deduce that
\begin{align*}
\sum_{w\in L(F_n')}\Omega_{F_n'}(w,t)&=\Omega_{T_n}(w_n,t)\prod _{i=1}^{n-1} \left(\Omega_{T_i}(w_{i}^{(1)},t)+\Omega_{T_i}(w_{i}^{(2)},t)+
\Omega_{T_i}(w_{i}^{(3)},t)+\Omega_{T_i}(w_{i}^{(4)},t)\right)\\
&=(t+1)^n(2t+1)^{n-1},
\end{align*}
as required.
 \qed

We now show that  the descent number over $P_n$ and $U_n$  is equidistributed with the
ascent number over $I'_{2n}$ and $I'_{2n-1}$, respectively.

\begin{thm}\label{theorem-3}
For $n\geq 1$, we have
\begin{equation}\label{main-m}
  \sum\limits_{\sigma\in P_n}x^{\mathrm{des}_{B} (\sigma)}
  =\sum\limits_{e\in I'_{2n}}x^{\mathrm{asc}(e)}
\end{equation}
and
\begin{equation}\label{main-n}
  \sum\limits_{\sigma\in U_n}x^{\mathrm{des}_{B} (\sigma)}
  =\sum\limits_{e\in I'_{2n-1}}x^{\mathrm{asc}(e)}.
\end{equation}
\end{thm}

To prove  Theorem \ref{theorem-3}, we need the following formulas \eqref{mmm} and \eqref{bn-2} of Savage and Schuster \cite{SavageM} for the generating function of the ascent number over $s$-inversion sequences of length $n$.
For a sequence $s=(s_1,s_2,\ldots)$
of positive integers, let $f_n^{(s)}(t)$ denote the number of sequences $(\lambda_1,\lambda_2,\ldots,\lambda_n)$ of nonnegative integers such that
\begin{equation}\label{bn-1}
0\leq \frac{\lambda_1}{s_1}\leq \frac{\lambda_2}{s_2}\leq\cdots \leq \frac{\lambda_n}{s_n}\leq t.
\end{equation}

\begin{thm}[\mdseries{\cite[Theorem 5]{SavageM}}]\label{bn-3}
Let $s=(s_1,s_2,\ldots)$ be a sequence of positive integers.
Then
\begin{equation}\label{mmm}
\sum_{t\geq 0} f_n^{(s)} (t)\,x^t = {1 \over (1-x)^{n+1}} \sum_{e} x^{\mathrm{asc}(e)},
\end{equation}
where $e$ ranges over   $s$-inversion sequences of length $n$.
\end{thm}

In particular, for  $s=(2,2,6,4,10,6,\ldots)$ and  \[s'=s/2=(1,1,3,2,5,3,\ldots),\]
Savage and Schuster \cite[Theorem 14]{SavageM}
showed   that
\begin{equation}\label{bn-2}
f_n^{(s')}(t)=(t+1)^{\lceil\frac{n}{2}\rceil}
\left(\frac{t+2}{2}\right)^{\lfloor\frac{n}{2}\rfloor}.
\end{equation}

We now turn to the proof of Theorem
\ref{theorem-3}.

\noindent
\textit{Proof of Theorem \ref{theorem-3}.}
It follows from \eqref{bn-2}  that
\[f_n^{(s)}(t)=f_n^{(s')}(2t)=\left(t+1\right)^{\lfloor\frac{n}{2}\rfloor}
(2t+1)^{\lceil\frac{n}{2}\rceil}.\]
By Theorem \ref{bn-3}, we obtain that
\begin{equation}\label{mainequation2}
\frac{\sum\limits_{e\in I'_{n}}x^{\mathrm{asc}(e)}}{(1-x)^{n+1}}=
\sum\limits_{t\geq0}\,(t+1)^{\lfloor\frac{n}{2}\rfloor}
(2t+1)^{\lceil\frac{n}{2}\rceil}\,x^t.
\end{equation}
To prove \eqref{main-m},
we replace  $n$ by $2n$ in   \eqref{mainequation2} to get
\begin{equation}\label{th361}
\frac{\sum\limits_{e\in I'_{2n}}x^{\mathrm{asc}(e)}}{(1-x)^{2n+1}}
=\sum\limits_{t\geq 0}\,(t+1)^n(2t+1)^n\,x^t.
\end{equation}
Combining the above relation and the formula \eqref{bkg} for $I_{2n}(x)$, we obtain
\[\sum\limits_{e\in I'_{2n}}x^{\mathrm{asc}(e)}=I_{2n}(x).\]
Thus \eqref{main-m} follows from the fact that  $I_{2n}(x)=P_n(x)$.

Next we prove \eqref{main-n}. Replacing $n$ by $2n-1$ in   \eqref{mainequation2},
we get
 \begin{equation*}
\frac{\sum\limits_{e\in I'_{2n-1}}x^{\mathrm{asc}(e)}}{(1-x)^{2n}}
=\sum\limits_{t\geq 0}\,(t+1)^{n-1}(2t+1)^{n}\,x^t.
\end{equation*}
Hence   \eqref{main-n} can be deduced from the following relation
\begin{equation}\label{tired-2}
\frac{\sum\limits_{\sigma\in U_n}x^{\mathrm{des}_B(\sigma)}}{(1-x)^{2n}}
=\sum\limits_{t\geq 0}\,(t+1)^{n-1}(2t+1)^{n}\,x^t.
\end{equation}
To prove \eqref{tired-2},  let
\[ U_n(x)=\sum_{\sigma\in U_n}x^{\mathrm{des}_B(\sigma)}.\]
We claim that $U_n(x)$ coincides with  the generating
function of the descent number over linear extensions of certain signed labeled
forests.
In the notation  $F_n', T_1,\ldots, T_n$
as defined in the proof of Theorem \ref{theorem-2},
we use  $\overline{L}(F_n')$ to denote  the set of signed labelings $w$ of $F_n'$
such that $w(T_n)=2n-1$ or $-(2n-1)$ and the induced labeling of $w$ on $F_{n-1}$ belongs to $L(F_{n-1})$.
Let
\[
G_n(x)=\sum_{w\in \overline{L}(F_n')}\,\sum_{\sigma\in\mathcal{L}(F_n',w)}x^{\mathrm{des}_B(\sigma)}.
\]
Again, using the argument in the proof of Theorem \ref{al}, one can construct a descent preserving bijection between the set
\[
\left\{\sigma\in \mathcal{L}(F_n',w)\,|\, w\in \overline{L}(F_n')\right\}
\]
and the set $U_n$. It then follows that
\begin{equation}\label{aj}
G_n(x)=U_n(x).
\end{equation}
Therefore,   \eqref{tired-2} is equivalent to
\begin{equation}\label{llk}
\frac{G_n(x)}{(1-x)^{2n}}
=\sum\limits_{t\geq 0}\,(t+1)^{n-1}(2t+1)^{n}\,x^t.
\end{equation}
Applying  Theorem \ref{pppt} to the set of signed labeled forests $(F_n',w)$ with $w\in \overline{L}(F_n')$, we find that
\[
\frac{G_n(x)}{(1-x)^{2n}}=\sum\limits_{t\geq 0}
\sum_{w\in \overline{L}(F_n')}\Omega_{F_n'}(w,t)\,
x^t.
\]
Hence \eqref{llk} can be deduced from the following relation
\begin{equation}\label{nnn}
\sum_{w\in \overline{L}(F_n')}\Omega_{F_n'}(w,t)
=(t+1)^{n-1}(2t+1)^{n}.
\end{equation}

To prove \eqref{nnn}, for $1\leq i\leq n-1$
and $1\leq j\leq 4$,
let $w_{i}^{(j)}$ be the signed labeling
of $T_i$ that is of type $j$, and let $w_{n}'$ and $w_n''$ be the signed labelings of $T_n$ such that $w_{n}'(T_n)=2n-1$ and  $w_{n}''(T_n)=-(2n-1)$.
As shown in the proof of Theorem \ref{new-2},
for $1\leq i\leq n-1$,
\[
 \Omega_{T_i}(w_{i}^{(1)},t)+\Omega_{T_i}(w_{i}^{(2)},t)+
\Omega_{T_i}(w_{i}^{(3)},t)+\Omega_{T_i}(w_{i}^{(4)},t)
=(t+1)(2t+1).
\]
Evidently, $\Omega_{T_n}(w_n',t)=t$ and $\Omega_{T_n}(w_n'',t)=t+1$.
Hence the sum on the left-hand side of \eqref{nnn} equals
\begin{align*}
&\left(\Omega_{T_n}(w_n',t)+\Omega_{T_n}(w_n'',t)\right)
\prod_{i=1}^{n-1}\left(\Omega_{T_i}(w_{i}^{(1)},t)+\Omega_{T_i}(w_{i}^{(2)},t)+
\Omega_{T_i}(w_{i}^{(3)},t)+\Omega_{T_i}(w_{i}^{(4)},t)\right)\\[3pt]
&\quad \quad=(t+1)^{n-1}(2t+1)^{n},
\end{align*}
as required.
\qed

\vspace{0.5cm}
 \noindent{\bf Acknowledgments.}  This work was supported by  the 973
Project and the National Science
Foundation of China.

\end{document}